\documentclass{llncs}

\usepackage[dvips]{graphicx}
\usepackage{amsmath,amssymb}
\usepackage{eclbkbox}
\usepackage{ascmac}

\newcommand{\diff}{\mathrm{d}}

\begin{document}
\title{Computational algebraic methods in efficient estimation}
%\title{Asymptotically Efficient Estimators for Algebraic Statistical Manifolds}

\author{Kei Kobayashi \and Henry P. Wynn}

\institute{The Institute of Statistical Mathematics  \email{kei@ism.ac.jp} \and London School of Economics \email{H.Wynn@lse.ac.uk}}

%\thanks{This work was supported by JSPS KAKENHI Grant 24700288 and
% EPSRC Grant EP/H007377/1}

\maketitle

\begin{abstract}
A strong link between information geometry and algebraic statistics is made by
investigating statistical manifolds which are algebraic varieties. In particular it it shown how
first and second order efficient estimators can be constructed, such as 
bias corrected Maximum Likelihood and more general estimators, and for which the estimating equations are 
purely algebraic. In addition it is shown how Gr\"obner basis technology, which is at the heart of 
algebraic statistics, can be used to reduce the degrees of the terms in the estimating equations.
This points the way to the feasible use, to find the estimators, of special methods for solving 
polynomial equations, such as homotopy continuation methods. 
Simple examples are given showing both equations and computations.
\end{abstract}

\section{Introduction}

Information geometry gives geometric insights and methods for studying the statistical efficiency of 
estimators, testing, prediction and model selection.
The field of algebraic statistics has proceeded somewhat separately but recently
a positive effort is being made to bring the two subjects together, notably \cite{gibilisco_etal2009}. 
This paper should be seen as part of this effort.  

A straightforward way of linking the two areas is to ask how far algebraic methods can be used when the
statistical manifolds of information geometry are algebraic, that is algebraic varieties or derived forms,
such as rational quotients. We call such models ``algebraic statistical models'' and will give formal definitions.

In the standard theory for non-singular statistical models, maximum likelihood estimators (MLEs) have 
first-order asymptotic efficiency and bias-corrected MLEs have second-order asymptotic efficiency. 
A short section covers briefly the basic theory of asymptotic efficiency using  differential geometry, necessary for our development.

We shall show that for some important algebraic models, the estimating equations of MLE type become polynomial and the degrees usually become very high if the model has a high-dimensional parameter space.
In this paper, asymptotically efficient algebraic estimators, a generalization of bias corrected MLE, are studied. By algebraic
estimators we mean estimators which are the solution of of algebraic equations. 
A main result is that for (algebraic) curved exponential family, there are second-order efficient estimators whose
polynomial degree is at most two.   
These are computed  by decreasing the degree of the estimating equations using Gr\"obner basis methods, the main tool of algebraic statistics. We supply  some the basic Gr\"obner theory in Appendix A. 
See \cite{pistone1996}.

The reduction of the degree saves computational costs dramatically when we use computational methods for solving the 
algebraic estimating equations. Here we use homotopy continuation methods of 
\cite{verschelde1999} \cite{li1997} to demonstrate this
effect for a few simple examples, for which we are able to carry out the Gr\"obner basis reduction.  
Appendix B discusses homotopy continuation methods.

Although, as mentioned, the links between computational algebraic methods and the theory of efficient estimators based on differential geometry are recent, two other areas of statistics, not covered here, exploit differential geometry methods. The first is tube theory.  The seminal paper by \cite{weyl1939} has been used to give exact confidence level values (size of tests), and bounds, for certain Gaussian simultaneous inference problems: \cite{naiman}, \cite{kuriki2002}.  This is very much related to the theory of up-crossings of Gaussian processes using expected Euler characteristic methods, see \cite{adler} and earlier papers. The second area is the use of the resolution of singularities (incidentally related to the tube theory) in which confidence levels are related to the dimension and the solid angle tangent of cones with apex at a singularity in parameters space \cite{drton}, \cite{watanabe}.
Moreover, the degree of estimating equations for MLE has been studied for some specific algebraic models, which are not necessarily singular \cite{drton_book}.
 In this paper
we cover the non-singular case, for rather more general estimators than MLE, 
and show that algebraic methods have a part to play.  

Most of the theories in the paper can be applied to a wider class of Multivariate Gaussian models with
some restrictions on their covariance matrices, for example models studied
in \cite{andersson-1998} \cite{gehrmann-2012}. Though the second-order efficient estimators proposed 
in the paper can be applied to them potentially, the cost for computing Gr\"obner basis
prevents their direct application.
Further innovation in the algebraic computation is required for real applications, which is a feature
of several other areas of algebraic statistics.

The next section gives some basic background in estimation and differential geometry for it.
Sections 3 and 4, which are the heart of the paper, give the algebraic developments and 
Section 5 gives some examples. 
Section 6 carries out some computation using homotopy continuation.

\section{Statistical manifolds and efficiency of estimators}

In this section, we introduce the standard setting of 
statistical estimation theory, via information geometry.
See \cite{amari1985} and \cite{amari-nagaoka2007} for details.
It is recognized that the ideas go back to at least
the work of Rao \cite{rao1945}, Efron \cite{efron1975} and Dawid \cite{dawid1977}.
The subject of information geometry was initiated by Amari and his collaborators
\cite{amari1982}, \cite{amari1983}.

Central to this family of ideas is that the rates of convergence of statistical estimators
and other test statistics depend on the metric and curvature of the parametric manifolds in
a neighborhood of the MLE or the null hypothesis.
In addition Amari realized the importance of two special models, the affine exponential model
and the affine mixture model, $e$ and $m$ frame respectively.
In this paper we concentrate on the exponential family model
but also look at curved subfamilies.
By extending the dimension of the parameter space of the exponential family,
we are able to cover some classes of mixture models.
The extension of the exponential model to infinite dimensions
is covered by\cite{pistone1995}.

\subsection{Exponential family and estimators}
\label{sec:exp-est}

\begin{it}A full exponential family \end{it} is a set of probability distributions
$\{\diff P(x|\theta) \mid \theta \in \Theta \}$ with a
parameter space $\Theta \subset \mathbb{R}^d$ such that
$$\diff P(x|\theta)=\exp(x_i \theta^i-\psi(\theta)) \diff\nu,$$
where $x \in \mathbb{R}^d$ is a variable representing a sufficient statistic
and $\nu$ is a carrier measure on $\mathbb{R}^d$.
Here $x_i \theta^i$ means $\sum_{i} x_i \theta^i$ (Einstein summation notation).

We call $\theta$ a natural parameter and 
$\eta=\eta(\theta):=E[x|\theta]$ an expectation parameter.
Denote $E= E(\Theta):=\{\eta(\theta) \mid \theta\in \Theta\}\subset \mathbb{R}^d$
as the corresponding expectation parameter space.
Note that the relation $\eta(\theta)=\nabla_\theta \psi(\theta)$ holds.
If the parameter space is restricted to a subset $\mathcal{V}_\Theta\subset \Theta$, 
we obtain a \begin{it}curved exponential family \end{it} 
$$\{\diff P(x|\theta) \mid \theta \in \mathcal{V}_\Theta \}.$$
The corresponding space of the expectation parameter is denoted by
$\mathcal{V}_E:=\{\eta(\theta) \mid \theta\in \mathcal{V}_\Theta\}\subset E$.
%We usually call such $\mathcal{V}_\Theta$, $\mathcal{V}_E$ and corresponding set of distributions
%endowed with the Fisher metric as model manifolds.

\begin{figure}
\begin{center}
\includegraphics[height=5cm]{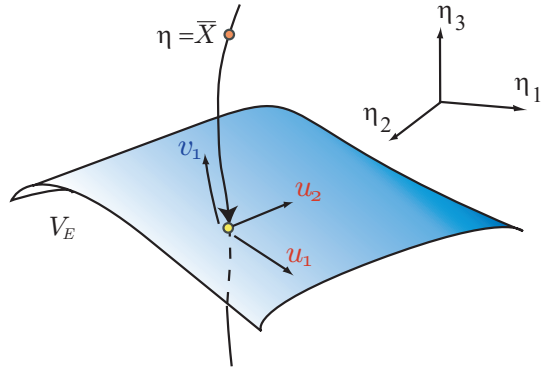}
\caption{A projection to the model manifold according to a local coordinate
defines an estimator.}
\label{fig:coordinates}
\end{center}
\end{figure}

Figure \ref{fig:coordinates} explains how to define an estimator by a local coordinate.
Let
$(u,v)\in \mathbb{R}^p \times \mathbb{R}^{d-p}$ 
with a dimension $p$ of $\mathcal{V}_\Theta$
be a local coordinate system around the true parameter $\theta^*$
and define $\mathcal{U} \subset \mathbb{R}^p$
such that $\{\theta(u,0)| u\in \mathcal{U}\}=\mathcal{V}_\Theta$.
For a full exponential model with $N$ samples obtained
by composing a map $(X^{(1)},\dots, X^{(N)}) \mapsto \theta(\eta)|_{\eta=\bar{X}}$
and a coordinate projection map $\theta(u,v)\mapsto u$,
we can define a (local) estimator $(X^{(1)},\dots, X^{(N)})\mapsto u$.
We define an estimator by $\eta(u,v)$ similarly.
Since $\bar{X}$ is a sufficient statistic of $\theta$ (and $\eta$)
in the full exponential family, every estimator can be computed by
$\bar{X}$ rather than the original data $\{X_i\}$.
Therefore in the rest of the paper, we write $X$ as shorthand for $\bar{X}$.

\subsection{Differential geometrical entities}
\label{subsec:geo-obj}
Let $w:=(u,v)$ and use indexes
$\{i,j,...\}$ for $\theta$ and $\eta$, $\{a,b,...\}$ for $u$, $\{\kappa,\lambda,...\}$ 
for $v$ and $\{\alpha,\beta,...\}$ for $w$.
The following  are used for expressing
conditions for asymptotic efficiency of estimators,
where Einstein notation is used.

\begin{itembox}[l]{Differential geometrical entities}
\begin{itemize}
\item $\eta_i(\theta) = \frac{\partial}{ \partial \theta^i}\psi(\theta)$,
\item Fisher metric $\underbar{G}=(g_{ij})$ w.r.t. $\theta$: $g_{ij}(\theta)=\frac{\partial^2 \psi(\theta)}{\partial \theta^i \partial \theta^j}$,
\item Fisher metric $\bar{G}=(g^{ij})$ w.r.t. $\eta$: $\bar{G}=\underbar{G}^{-1}$,
\item Jacobian: $B_{i\alpha}(\theta):=\frac{\partial \eta_i(w)}{\partial w^\alpha}$,
\item e-connection: $\Gamma^{(e)}_{\alpha\beta,\gamma}=(\frac{\partial^2}{\partial w^\alpha \partial w^\beta} \theta^i(w)) (\frac{\partial}{\partial w^\gamma} \eta_i(w)) $,
\item m-connection: $\Gamma^{(m)}_{\alpha\beta,\gamma}=(\frac{\partial^2}{\partial w^\alpha \partial w^\beta} \eta_i(w)) (\frac{\partial}{\partial w^\gamma} \theta^i(w)) $,
\end{itemize}
\end{itembox}

\subsection{Asymptotic statistical inference theory} 
Under some regularity conditions
on the carrier measure $\nu$, potential function $\psi$ 
and the manifolds $\mathcal{V}_\Theta$ or $\mathcal{V}_E$,
the asymptotic theory below is available. 
These conditions are for guaranteeing 
the finiteness of the moments and the commuting of the expectation and the partial derivative
$\frac{\partial}{\partial \theta} E_\theta[f]= E_\theta[\frac{\partial f}{\partial\theta}]$.
For more details of the required regularity conditions, see Section 2.1 of \cite{amari1985}.

\begin{enumerate}
\item If $\hat{u}$ is a consistent estimator (i.e. $P(\|\hat{u}-u\|>\epsilon)\rightarrow 0$ as $N\rightarrow \infty$ for any $\epsilon>0$),
the squared error matrix of $\hat{u}$ is
$$E_u[(\hat{u}-u)(\hat{u}-u)^\top]=E_u[(\hat{u}^a-u^a)(\hat{u}^b-u^b)]=N^{-1}[g_{ab}-g_{a\kappa}g^{\kappa \lambda}g_{b\lambda}]^{-1}+O(N^{-2}).$$
Here $[\cdot]^{-1}$ means the matrix inverse.
Thus, if $g_{a\kappa}=0$ for all $a$ and $\kappa$, the main term in the r.h.s. becomes minimum. We call such an estimator as a \begin{it} 1-st order efficient \end{it} estimator.
\item
The bias term becomes 
$$E_u[\hat{u}^a-u^a]=(2N)^{-1}b^a(u)+ O(N^{-2})$$
 for each $a$ where $b^a(u):=\Gamma^{(m)}{}_{cd}^a(u) g^{cd}(u)$.
Then, the \begin{it} bias corrected estimator \end{it}
$\check{u}^a:=\hat{u}^a-b^a(\hat{u})$ satisfies $E_u[\check{u}^a-u^a]=O(N^{-2})$.
\item
Assume $g_{a\kappa}=0$ for all $a$ and $\kappa$, then the square error matrix is represented by
\begin{equation*}
E_u[(\check{u}^a-u^a)(\check{u}^b-u^b)]=\frac{1}{N}g^{ab}+\frac{1}{2N^2}
(\mathop{{\Gamma}^{2ab}_M}^{(m)}+2\mathop{{H}^{2ab}_M}^{(e)}+\mathop{{H}^{2ab}_A}^{(m)})+o(N^{-2}).
\end{equation*}
See Theorem 5.3 of \cite{amari1985} and Theorem 4.4 of \cite{amari-nagaoka2007}
for the definition of the terms in the r.h.s. 
Of the four dominating terms in the r.h.s., only
\begin{equation*}
\mathop{{H}^{2ab}_A}^{(m)}:=g^{\kappa\mu}g^{\lambda\nu}
H^{(m)}{}_{\kappa\lambda}^a H^{(m)}{}_{\mu\nu}^b
%\label{eq:H2ab}
\end{equation*}
depends on the selection of the estimator.

Here $H^{(m)}{}_{\kappa\lambda}^a$ is an embedding curvature and 
equal to $\Gamma^{(m)}{}_{\kappa\lambda}^a$ when $g_{a\kappa}=0$ for every $a$ and $\kappa$.
Since $\mathop{{H}^{2ab}_A}^{(m)}$ is the square of $\Gamma^{(m)}{}_{\kappa\lambda}^a$, 
the square error matrix attains the minimum in the sense of positive definiteness if and only if
\begin{equation}
\label{eq:2nd-eff}
\left.\Gamma^{(m)}{}_{\kappa\lambda,a}(w)\right|_{v=0}=\left.\left(\frac{\partial^2}{\partial v^\kappa \partial v^\lambda}
\eta_i(w)\right) \left(\frac{\partial}{\partial u^a} \theta^i(w)\right)\right|_{v=0}=0.
\end{equation}
Therefore the elimination of the m-connection (\ref{eq:2nd-eff}) implies \begin{it} second-order efficiency
\end{it} of the estimator after a bias correction, i.e.
it becomes optimal among the bias-corrected first-order efficient estimators
up to $O(N^{-2})$.
\end{enumerate}

\section{Algebraic models and efficiency of algebraic estimators}
\label{sec:estimation}

This section studies asymptotic efficiency for
statistical models and estimators which are defined algebraically.
Many models in statistics are defined algebraically.
Perhaps most well known are polynomial regression models and 
algebraic conditions on probability models such as independence
and conditional independence. 
Recently there has been considerable interest
in marginal models \cite{bergsma2009} which are typically
linear restrictions on raw probabilities.
In time series autoregressive models expressed by linear transfer functions
induce algebraic restrictions on covariance matrices.
Our desire is to have a definition of algebraic statistical model
which can be expressed from within the curved exponential family
framework but is sufficiently broad to cover cases
such as those just mentioned.
Our solution is to allow algebraic conditions in the natural parameter $\theta$,
mean parameter $\eta$ or both.
The second way in which algebra enters is in the form of the estimator.

\subsection{Algebraic curved exponential family}
We say a curved exponential family is \begin{it} algebraic \end{it} if
the following two conditions are satisfied.

\begin{itemize}

\item[(C1)]
$\mathcal{V}_\Theta$ or $\mathcal{V}_E$ is represented by a real algebraic variety,
i.e. $\mathcal{V}_\Theta:=\mathcal{V}(f_1,\dots,f_k)=\{\theta\in \mathbb{R}^d | f_1(\theta)=\dots=f_k(\theta)=0\}$ or similarly $\mathcal{V}_E :=\mathcal{V}(g_1,\dots,g_k)$
for
$f_i \in \mathbb{R}[\theta^1,\dots,\theta^d]$ and
$g_i \in \mathbb{R}[\eta_1,\dots,\eta_d]$.

\item[(C2)]
$\theta \mapsto \eta(\theta)$ or $\eta \mapsto \theta(\eta)$ is represented by some algebraic equations,
i.e. there are $h_1,\dots,h_k \in \mathbb{R}[\theta,\eta]$ such that
locally in $\mathcal{V}_\Theta \times \mathcal{V}_E$,
$h_i(\theta,\eta)=0$ iff $\eta(\theta)=\eta$ or $\theta(\eta)=\theta$.
\end{itemize}

Here $\mathbb{R}[\theta^1,\dots,\theta^d]$ means a polynomial of $\theta^1,\dots,\theta^d$ over
the real number field $\mathbb{R}$ and
$\mathbb{R}[\theta,\eta]$ means $\mathbb{R}[\theta^1,\dots,\theta^d,\eta_1,\dots,\eta_d]$.
The integer $k$, the size of the generators, is not necessarily equal to $d-p$ but
we assume $\mathcal{V}_\Theta$ (or $\mathcal{V}_E$) has dimension $p$ around
the true parameter.
Note that if $\psi(\theta)$ is a rational form or the logarithm of a rational form,
(C2) is satisfied.

\subsection{Algebraic estimators}
\label{sec:alg-est}
The parameter set $\mathcal{V}_\Theta$ (or $\mathcal{V}_E$) 
is sometimes singular for algebraic models.  
But throughout the following analysis, we assume non-singularity around the true parameter $\theta^*\in \mathcal{V}_\Theta$ (or $\eta^*\in \mathcal{V}_E$ respectively) .

Following the discussion at the end of Section \ref{sec:exp-est}.
We call $\theta(u,v)$ or $\eta(u,v)$ an \begin{it}algebraic estimator\end{it} if
\begin{itemize}
\item[(C3)]
 $w \mapsto \eta(w)$ or $w \mapsto \theta(w)$ is represented algebraically.
\end{itemize}
We remark that the MLE for an algebraic curved exponential family is an algebraic estimator.

If conditions (C1), (C2) and (C3) hold, then all of the geometrical entities
in section \ref{subsec:geo-obj} are characterized by special polynomial equations.
Furthermore, if $\psi(\theta)\in \mathbb{R}(\theta) \cup \log \mathbb{R}(\theta)$
and $\theta(w)\in \mathbb{R}(w) \cup \log \mathbb{R}(w)$,
then the geometrical objects have the additional property of being rational.

\subsection{Second-order efficient algebraic estimators, vector version}

Consider an algebraic estimator $\eta(u,v)\in \mathbb{R}[u,v]^d$
satisfying the following vector equation:
\begin{equation}
\label{eq:2nd-eff-explicit}
X=\eta(u,0)+\sum_{i=p+1}^{d} v_{i-p} e_i(u)+ c\cdot \sum_{j=1}^p f_j(u,v) e_j(u)
\end{equation}
where, for each $u$, $\{e_j(u); j=1,\dots,p\}\cup \{e_i(u); i=p+1,\dots,d\}$ is a complete basis of $\mathbb{R}^d$
such that
$\langle e_j(u),(\bigtriangledown_u \eta)\rangle_g=0$
and $f_j(u,v)\in \mathbb{R}[u][v]_{\geq 3}$, namely a polynomial whose degree in $v$ is at least 3
with coefficients polynomial in $u$, for $j=1,\dots,p$. 
Remember we use a notation $X=\bar{X}=\frac{1}{N}\sum_i X_i$.
The constant $c$ is to control the perturbation (see below).
% When the constant $c=0$, 
% we obtain the MLE.
% Otherwise $c\in \mathbb{R}$ may be freely chosen
% and represent a one-dimensional perturbation from MLE.

A straightforward computation of the $m$-connection in (\ref{eq:2nd-eff}) 
at $v=0$ for
$$\eta(w)=\eta(u,0)+\sum_{i=p+1}^{d} v_{i-p} e_i(u)+ c\cdot \sum_{j=1}^p f_j(u,v) e_j(u)$$
shows it to be zero. This gives
\begin{theorem}
\label{thm:2to1}
Vector equation (\ref{eq:2nd-eff-explicit}) satisfies the second-order efficiency 
(\ref{eq:2nd-eff}).
\end{theorem}

We call (\ref{eq:2nd-eff-explicit}) \begin{it}a vector version of a second-order efficient
estimator\end{it}.
Note that if $c=0$, (\ref{eq:2nd-eff-explicit}) gives an
estimating equation for the MLE. Thus the last term in (\ref{eq:2nd-eff-explicit})
can be recognized as a perturbation from the MLE.

\begin{figure}
\begin{center}
\includegraphics[height=5cm]{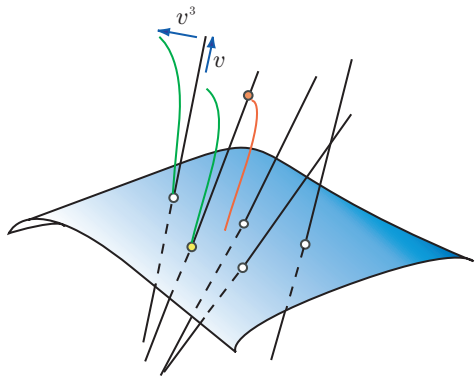}
%\pgfuseimage{real_regression_nobias3.eps}
\caption{Image of the vector version of the second-order efficient estimators}
\label{fig:image-vec}
\end{center}
\end{figure}

Figure \ref{fig:image-vec} is a rough sketch of the second-order efficient estimators.
Here the model is embedded in an $m$-affine space.
Given a sample (red point), the MLE is an orthogonal projection (yellow point) to the model
with respect to the Fisher metric. But a second-order efficient estimator maps
the sample to the model along a ``cubically'' curved manifold (red curve).

\subsection{Second-order efficient algebraic estimators, algebraic version}
\label{subsec:alg-form}
%The solution of (\ref{eq:2nd-eff-explicit}) is given by a set of equations as follows:

Another class of second-order efficient algebraic estimators we call the
\begin{it} algebraic version\end{it},
which is defined by the following simultaneous polynomial equations with $\eta_u=\eta(u,0)$.
\begin{align}
(X-\eta_u)^\top \tilde{e}_j(u,\eta_u)&+c\cdot h_j(X,u,\eta_u,X-\eta_u)=0 
\mbox{~for~} j=1,\dots,p \label{eq:2nd-eff-implicit}
\end{align}
where
$\{\tilde{e}_j(u,\eta_u)\in \mathbb{R}[u,\eta_u]^d; j=1,\dots,p\} $ span $((\nabla_u \eta(u,0))^{\perp_{\bar{G}}})^{\perp_E}$
 for every $u$ and
$h_j(X,u,\eta_u,t) \in \mathbb{R}[X,u,\eta_u][t]_3$ ($\mbox{degree}=3$ in $t$) for $j=1,\dots,p$.
The constant $c$ is to control the perturbation.
The notation $\bar{G}$ represents the Fisher metric on the full-exponential family with respect to $\eta$.
The notation $(\nabla_u \eta(u,0))^{\perp_{\bar{G}}}$ means the subspace orthogonal to $span(\partial_a \eta(u,0))_{a=1}^p$
with respect to $\bar{G}$ and $(\cdot)^{\perp_E}$ means the orthogonal complement in the sense of
Euclidean vector space.
Here, the term ``degree'' of a polynomial means the maximum degree of its terms.
Note that the case
$(X-\eta_u)^\top \tilde{e}_j(u,\eta_u)=0$ for $j=1,\dots,p$ gives a special set of the estimating equations of the MLE.

\begin{theorem}
\label{thm:2to3}
An estimator defined by a vector version (\ref{eq:2nd-eff-explicit})  of the second-order efficient estimators is also represented by an algebraic version (\ref{eq:2nd-eff-implicit})
where $h_j(X,u,\eta_u,t)=\tilde{f}_j(u,(\tilde{e}_i^\top t)_{i=1}^p,(\tilde{e}_i^\top (X-\eta_u))_{i=1}^p)$
with a function $\tilde{f}_j(u,v,\tilde{v})\in \mathbb{R}[u,\tilde{v}][v]_3$ such that
$\tilde{f}(u,v,v)=f(u,v)$. 
\end{theorem}
\proof
Take the Euclidean inner product of both sides of (\ref{eq:2nd-eff-explicit})  with each $\tilde{e}_j$ which is a vector Euclidean orthogonal to the subspace $span(\{e_i |i\neq j\})$ and obtain a system of polynomial equations. 
By eliminating variables $v$ from the polynomial equations, an algebraic version is obtained.
\qed
\vspace{0.5cm}

% This can be proved by (i) taking the Euclidean inner product of both sides of (\ref{eq:2nd-eff-explicit})  with each $\tilde{e}_j$ which is an vector Euclidean orthogonal to the subspace $span(\{e_i |i\neq j\})$ and (ii) removing variables $v$.

\begin{theorem}
\label{thm:3to1}
Every algebraic equation
(\ref{eq:2nd-eff-implicit}) gives a second-order efficient estimator (\ref{eq:2nd-eff}).
\end{theorem}
\proof
Writing  $X=\eta(u,v)$ in (\ref{eq:2nd-eff-implicit}), we obtain
$$(\eta(u,v)-\eta(u,0))^\top \tilde{e}_j(u) + c \cdot h_j(\eta(u,v),u,\eta(u,0),\eta(u,v)-\eta(u,0))=0.$$
Partially differentiate this by $v$ twice, we obtain
$$\left.\left(\frac{\partial^2 \eta(u,v)}{\partial v^\lambda \partial v^\kappa}\right)^\top \tilde{e}_j(u)\right|_{v=0} = 0,$$
since each term of $h_j(\eta(u,v),u,\eta(u,0),\eta(u,v)-\eta(u,0))$ has degree more than 3 
in its third component $(\eta_i(u,v)-\eta_i(u,0))_{i=1}^{d}$ and
$\left.\eta(u,v)-\eta(u,0)\right|_{v=0}=0$.
Since ${\rm span}\{\tilde{e}_j(u); j=1,\dots,p\}=((\nabla_u \eta(u,0))^{\perp_{\bar{G}}})^{\perp_E}={\rm span}\{\bar{G} \partial_{u_a} \eta
;a=1,\dots,p\}$, we obtain
$$\left.\Gamma_{\kappa\lambda a}^{(m)} \right|_{v=0} =
\left.\frac{\partial^2 \eta_i}{\partial v^\lambda \partial v^\kappa}
g^{ij} \frac{\partial \eta_j}{\partial u^a} \right|_{v=0} = 0.$$
This implies the estimator is second-order efficient. \qed
\vspace{0.5cm}

By Theorems \ref{thm:2to1}, \ref{thm:2to3} and \ref{thm:3to1},
the relationship between the three forms of the second-order efficient algebraic
estimators is summarized as
$$(\ref{eq:2nd-eff}) \Leftarrow (\ref{eq:2nd-eff-explicit}) \Rightarrow (\ref{eq:2nd-eff-implicit}) \Rightarrow (\ref{eq:2nd-eff}).$$
Furthermore, if we assume the estimator has a form $\eta\in \mathbb{R}(u)[v]$, that is 
a polynomial in $v$ with coefficients rational in $u$,
every first-order efficient estimator satisfying (\ref{eq:2nd-eff}) can be written
in a form (\ref{eq:2nd-eff-explicit}) after resetting coordinates $v$ for the estimating manifold.
In this sense, we can say $(\ref{eq:2nd-eff}) \Rightarrow (\ref{eq:2nd-eff-explicit})$ and
the following corollary holds.
\begin{corollary}
If $\eta\in \mathbb{R}(u)[v]$,  the forms (1), (2) and (3) are equivalent.
\end{corollary}

\subsection{Properties of the estimators}
The following theorem is a straightforward extension of 
the local existence of MLE.
That is to say, the existence for sufficiently large sample size.
The regularity conditions are essentially the same as for the MLE
but with an additional condition referring to the control constant $c$.
\begin{proposition}[Existence and uniqueness of the estimate]
Assume that the Fisher matrix is non-degenerate around $\eta(u^*)\in \mathcal{V}_E$.
Then the estimate given by (\ref{eq:2nd-eff-explicit}) or (\ref{eq:2nd-eff-implicit}) locally uniquely exists 
for small $c$, i.e.
there is a neighborhood $G(u^*)\subset \mathbb{R}^d$ of
$\eta(u^*)$ and $\delta>0$ such that for every fixed $X \in G(u^*)$ and $-\delta<c<\delta$,
a unique estimate exists.
\end{proposition}
\proof
Under the condition of the theorem, the MLE always exists locally.
Furthermore, because of the nonsingular Fisher matrix, the MLE is locally bijective
(by the implicit representation theorem).
Thus $(u_1,\dots,u_p)\mapsto (g_1(x-\eta_u),\dots,g_p(x-\eta_u))$ 
for $g_j(x-\eta_u):=(X-\eta_u)^\top \tilde{e}_j(u,\eta_u)$
in (\ref{eq:2nd-eff-implicit}) is locally bijective. Since $\{g_i\}$ and $\{h_i\}$
are continuous, we can select $\delta>0$ for (\ref{eq:2nd-eff-implicit}) to be
locally bijective for every  $-\delta<c<\delta$. \qed

\subsection{Summary of estimator construction}

We summarize how to define a second-order efficient algebraic estimator
(vector version) and how to compute an algebraic version from it. 
\vspace{0.5cm}

\begin{breakbox}
{\bf Input}:
\vspace{-0.3cm}
\begin{itemize}
%$\psi\in \mathbb{R}(\eta) \cup \log \mathbb{R}(\eta)$, 
\item[$\cdot$] a potential function $\psi$ satisfying (C2), 
\item[$\cdot$]
polynomial equations of $\eta,u$ and $v$ satisfying (C3),
\item[$\cdot$]
$m_1,\dots,m_{d-p}\in \mathbb{R}[\eta]$ such that
$\mathcal{V}_E=V(m_1,\dots,m_{d-p})$ gives the model,
\item[$\cdot$]
$f_j \in\mathbb{R}[u][v]_{\geq 3}$ and $c\in \mathbb{R}$ for a vector version
\end{itemize}

\begin{itemize}
\item[]
Step\;1.
Compute $\psi$ and $\theta(\eta)$, $G(\eta)$, ($\Gamma^{(m)}(\eta)$ for bias correction)

\item[]
Step\;2.
Compute $f_{ai} \in \mathbb{R}[\eta][\xi_{11},\dots,\xi_{pd}]_1$ s.t. 
$f_{aj}(\xi_{11},\dots,\xi_{pd}):=\partial_{u^a} m_j$ for $
\xi_{b i} :=\partial_{u^b}\eta_i.$

\item[]
Step\;3.
Find $e_{p+1},\dots,e_{d}\in (\nabla_u \eta)^{\perp_{\bar{G}}}$ by eliminating $\{\xi_{aj}\}$ from
$\langle e_i,\partial_{u^a} \eta\rangle_{\bar{G}}=e_{ik}(\eta) g^{kj}(\eta) \xi_{aj}=0$
and $f_{aj}(\xi_{11},\dots,\xi_{pd})=0$.
%Syzygies may be used for finding $\{e_i(\eta)\}$.

\item[]
Step\;4.
Select $e_{1},\dots,e_p \in \mathbb{R}[\eta]$ s.t. $e_1(\eta),\dots,e_d(\eta)$ are
linearly independent.

\item[]
Step\;5.
Eliminate $v$ from
\begin{center}$X=\eta(u,0)+\sum_{i=p+1}^{d} v_{i-p} e_i+ c\cdot \sum_{j=1}^p f_j(u,v) e_j$\end{center}
and compute $(X-\eta)^\top \tilde{e}_j$ and $h\in (\mathbb{R}[\eta][X-\eta]_3)^p$,
given by Theorem \ref{thm:2to3}.
\end{itemize}

{\bf Output(Vector version)}: 
\begin{center}
$X=\eta(u,0)+\sum_{i=p+1}^{d} v_{i-p} e_i(\eta)+ c\cdot \sum_{j=1}^p f_j(u,v) e_j(\eta)$.
\end{center}
{\bf Output(Algebraic version)}: 
\begin{center}
$(X-\eta)^\top \tilde{e}+c \cdot h(X-\eta)=0$.
\end{center}
\end{breakbox}

\subsection{Reduction of the degree of the estimating equations}

As we noted in section \ref{subsec:alg-form}, if we set $h_j=0$ for all $j$,
the estimator becomes the MLE. In this sense, $c h_j$ can be recognized as a perturbation
from the likelihood equations. If we select each
$h_j(X,u,\eta_u,t) \in \mathbb{R}[X,u,\eta_u][t]_3$ tactically, we can reduce the degree of the
polynomial estimating equation.
For algebraic background, the reader refers to Appendix A.

Here, we assume $u\in \mathbb{R}[\eta_u]$. For example, we can set $u_i=\eta_i$.
Then $\tilde{e}_j(u,\eta_u)$ is a function of $\eta_u$, so
we write it as $\tilde{e}_j(\eta)$.
Define an ideal $\mathcal{I}_3$ of $\mathbb{R}[X,\eta]$ as
$$\mathcal{I}_3:=\langle\{(X_i-\eta_i)(X_j-\eta_j)(X_k-\eta_k)\mid 1\leq i,j,k \leq d \}\rangle.$$

Select a monomial order $\prec$ and set $\eta_1\succ\dots\succ\eta_d \succ X_1\succ \dots \succ X_d$.
Let $G_{\prec}=\{g_1,\dots,g_m\}$ be a Gr\"{o}bner basis of $\mathcal{I}_3$ with respect to $\prec$.
Then the remainder (normal form) $r_j$ of $(X-\eta)^\top \tilde{e}_j(\eta)$,
the first term of the l.h.s. of (\ref{eq:2nd-eff-implicit}),
with respect to $G_{\prec}$,  is uniquely determined for each $j$.

\begin{theorem}
If the monomial order $\prec$ is the pure lexicographic,
\begin{enumerate}
\item $r_j$ for $j=1,\dots, p$ has degree at most 2 with respect to $\eta$, and 
\item $r_j=0$ for $j=1,\dots, p$ are the estimating equations for a second-order efficient estimator.
\end{enumerate}
\end{theorem}
\proof
Assume $r_j$ has a monomial term whose degree is more than 2 with respect to $\eta$ and
represent the term as $\eta_a \eta_b \eta_c q(\eta,X)$ with a polynomial $q\in \mathbb{R}(\eta,X)$
and a combination of indices $a,b,c$.
Then $\{\eta_a \eta_b \eta_c + (X_a-\eta_a)(X_a-\eta_a)(X_a-\eta_a)\}q(\eta,X)$ has
a smaller polynomial order than $\eta_a \eta_b \eta_c q(\eta,X)$ since $\prec$ is
pure lexicographic satisfying $\eta_1\succ\dots\succ\eta_d \succ X_1\succ \dots \succ X_d$.
Therefore by subtracting  $(X_a-\eta_a)(X_a-\eta_a)(X_a-\eta_a)\}q(\eta,X)\in \mathcal{I}_3$ from $r_j$,
the polynomial degree decreases. 
This contradicts the fact $r_j$ is the normal form so each $r_j$ has degree at most 2.

Furthermore each polynomial in $\mathcal{I}_3$ is in $\mathbb{R}[X,u,\eta_u][X-\eta]_3$ and
therefore by taking the normal form, the condition for the algebraic version (\ref{eq:2nd-eff-implicit})
of second-order efficiency still holds.\qed

The reduction of the degree is important when we use algebraic algorithms such as 
homotopy continuation methods \cite{lee_etal2008} to solve simultaneous polynomial equations
since computational cost depends highly on the degree of the polynomials.

\section{First-order efficiency}
\label{sec:1st-eff}

It is not surprising that, for first-order efficiency, almost the same arguments hold
as for second-order efficiency.

By Theorem 5.2 of \cite{amari1985}, a consistent estimator is 1st-order efficient if and only if
\begin{equation}
\label{eq:1st-eff}
g_{\kappa a}=0.
\end{equation}

Consider an algebraic estimator $\eta(u,v)\in \mathbb{R}[u,v]^d$
satisfying the following vector equation:
\begin{equation}
\label{eq:1st-eff-explicit}
X=\eta(u,0)+\sum_{i=p+1}^{d} v_{i-p} e_i(u)+ c\cdot \sum_{j=1}^p f_j(u,v) e_j(u)
\end{equation}
where, for each $u$, $\{e_j(u); j=1,\dots,p\}\cup \{e_i(u); i=p+1,\dots,d\}$ is a complete basis of $\mathbb{R}^d$ s.t.
$\langle e_j(u),(\bigtriangledown_u \eta)\rangle_g=0$
and $f_j(u,v)\in \mathbb{R}[u][v]_{\geq 2}$, a polynomial whose degree of $v$ is
at least 2, for $j=1,\dots,p$. 
Similarly, $c\in \mathbb{R}$ is a constant for perturbation.  
Here, the only difference between  (\ref{eq:2nd-eff-explicit}) for the second-order efficiency
and (\ref{eq:1st-eff-explicit}) for the first-order efficiency is the degree of the $f_j(u,v)$
with respect to $v$. 

The algebraic version of the first-order efficient algebraic estimator
is defined by the following simultaneous polynomial equalities with $\eta_u=\eta(u,0)$.
\begin{align}
(X-\eta_u)^\top \tilde{e}_j(u,\eta_u)&+c\cdot h_j(X,u,\eta_u,X-\eta_u)=0 
\mbox{~for~} j=1,\dots,p \label{eq:1st-eff-implicit}
\end{align}
where
$\{\tilde{e}_j(u,\eta_u)\in \mathbb{R}[u,\eta_u]^d; j=1,\dots,p\} $ span $((\nabla_u \eta(u,0))^{\perp_{\bar{G}}})^{\perp_E}$
 for every $u$ and
$h_j(X,u,\eta_u,t) \in \mathbb{R}[X,u,\eta_u][t]_2$ ($\mbox{degree}=2$ w.r.t. $t$) for $j=1,\dots,p$.
Here, the only difference between  (\ref{eq:2nd-eff-implicit}) for the second-order efficiency
and (\ref{eq:1st-eff-implicit}) for the first-order efficiency is the degree of the $h_j(X,u,\eta_u,t)$
with respect to $t$. 

Then the relation between the three different forms of first-order efficiency can be proved in the same way manner as for Theorem \ref{thm:2to1}, \ref{thm:2to3} and \ref{thm:3to1}.
\begin{theorem}
(i) Vector version (\ref{eq:1st-eff-explicit}) satisfies the first-order efficiency.\\
(ii) An estimator defined by a vector version (\ref{eq:1st-eff-explicit}) of the first-order efficient estimators is also represented by an algebraic version (\ref{eq:1st-eff-implicit}).\\
(iii) Every algebraic version (\ref{eq:1st-eff-implicit}) gives a first-order efficient estimator.
\end{theorem}
The relationship between the three forms of the first-order efficient algebraic
estimators is summarized as
$(\ref{eq:1st-eff}) \Leftarrow (\ref{eq:1st-eff-explicit}) \Rightarrow (\ref{eq:1st-eff-implicit}) \Rightarrow (\ref{eq:1st-eff})$.
Furthermore, if we assume the estimator has a form $\eta\in \mathbb{R}(u)[v]$, 
the forms (\ref{eq:1st-eff}), (\ref{eq:1st-eff-explicit}) and (\ref{eq:1st-eff-implicit}) are equivalent.

Let $\mathcal{R}:=\mathbb{Z}[X,\eta]$ and define
$$\mathcal{I}_2:=\langle\{(X_i-\eta_i)(X_j-\eta_j)\mid 1\leq i,j \leq d \}\rangle$$
as an ideal of $\mathcal{R}$.
In a similar manner, let $\prec$ be a monomial order such that $\eta_1\succ\dots\succ\eta_d \succ X_1\succ \dots \succ X_d$.
Let $G_{\prec}=\{g_1,\dots,g_m\}$ be a Gr\"{o}bner basis of $\mathcal{I}_2$ with respect to $\prec$.
The properties of the normal form $r_i$ of $(X-\eta(u,0))^\top \tilde{e}_i(u)$ with respect to $G_{\prec}$ 
are then covered by the following:

\begin{theorem}
If the monomial order $\prec$ is the pure lexicographic,\\
(i) $r_i$ for $i=1,\dots, d$ has degree at most 1 with respect to $\eta$, and \\
(ii) $r_i=0$ for $i=1,\dots, d$ are the estimating equations for a first-order efficient estimator.
\end{theorem}

\section{Examples}
\label{sec:examples}
In this section, we show how to use the algebraic computation to design asymptotically efficient estimators
for two simple examples. The examples satisfy the algebraic conditions
(C1), (C2) and (C3) so it is verified that necessary geometric entities have an algebraic form
as mentioned in Section \ref{sec:alg-est}.

\subsection{Example: Periodic Gaussian Model}
The following periodic Gaussian model shows how to compute 
second-order efficients estimators and their biases.

\begin{itemize}
\item[\textbullet] Statistical Model:\\

$X\sim N(\mu, \Sigma(a))$
with
$\mu=
\begin{bmatrix}
0 \\
0 \\
0 \\
0
\end{bmatrix}$
and 
$\Sigma(a)=
\begin{bmatrix}
1 & a & a^2 & a \\
a & 1 & a & a^2 \\
a^2 & a & 1 & a \\
a & a^2 & a & 1
\end{bmatrix}
$ for $0\leq a < 1$.
\vspace{0.2cm}

Here, the dimension of the full exponential family and the curved exponential 
family are $d=3$ and $p=1$, respectively.

\item[\textbullet] Curved exponential family:
$$\log f(x|\theta) = 2 \left(x_1 x_2 + x_2 x_3 + x_3 x_4 + x_4 x_1\right)
\theta_{{2}}+ 2\left(x_{{3}}x_{{1}}+x_{{4}}x_{{2}}
 \right) \theta_{{3}}-\psi(\theta),$$
 
\item[\textbullet] Potential function:
$$\psi(\theta)=-1/2\,\log  ( {\theta_{{1}}}^{4}-4\,{\theta_{{1}}}^{2}{\theta_{{2}
}}^{2}+8\,\theta_{{1}}{\theta_{{2}}}^{2}\theta_{{3}}-2\,{\theta_{{1}}}
^{2}{\theta_{{3}}}^{2}-4\,{\theta_{{2}}}^{2}{\theta_{{3}}}^{2}+{\theta
_{{3}}}^{4} ) +2\,\log  ( 2\,\pi ),$$
\item[\textbullet] Natural parameter:
$$\theta(a)=\left[\frac{1}{1-2a^2+4a^4},-\frac{a}{1-2a^2+4a^4}, \frac{a^2}{1-2a^2+4a^4}\right]^\top,$$
\item[\textbullet] Expectation parameter:~~
$\eta(a)=[-2,-4a,-2a^2]^\top,$
\item[\textbullet] Fisher metric with respect to $\eta$:
$$(g^{ij})=\left[ {\small \begin {array}{ccc} 2\,{a}^{4}+4\,{a}^{2}+2&8\,a \left( 1+{a}^
{2} \right) &8\,{a}^{2}\\8\,a \left( 1+{a}^{2}
 \right) &4+24\,{a}^{2}+4\,{a}^{4}&8\,a \left( 1+{a}^{2} \right) 
\\ 8\,{a}^{2}&8\,a \left( 1+{a}^{2} \right) &2\,{a}^
{4}+4\,{a}^{2}+2\end {array} }\right] ,$$
\item[\textbullet] A set of vectors $e_i\in \mathbb{R}^3$:
$$e_0(a):=[0,-1,a]^\top \in \partial_a \eta(a),$$
$$e_1(a):=[3a^2+1,4a,0]^\top,~ e_2(a):=[-a^2-1,0,2]^\top \in (\partial_a \eta(a))^{\perp_{\bar{G}}}.$$

\item[\textbullet] A vector version of the second-order efficient estimator is, for example,
$$x-\eta+v_1\cdot e_1 + v_2 \cdot e_2 +c \cdot v_1^3\cdot e_0=0.$$

\item[\textbullet] A corresponding algebraic version of the second-order efficient estimator:
by eliminating $v_1$ and $v_2$, we get
$g(a)+c\cdot h(a)=0$ where
$$g(a):=8( a-1) ^{2} ( a+1) ^{2} ( 1+2{a}^{2}) ^{2}
( 4{a}^{5}-8{a}^{3}+2{a}^{3}{\it x_3}-3{
\it x_2}{a}^{2}+4a+4a{\it x_1}+2a{\it x_3}-{\it x_2})$$
$$\mbox{ and } h(a):=( 2{a}^{4}+{a}^{3}{\it x_2}-{a}^{2}{\it x_3}+2{a}^{2}+a{\it x_2
}-2{\it x_1}-{\it x_3}-4) ^{3}.$$

\item[\textbullet] An estimating equation for MLE:
$$4{a}^{5}-8{a}^{3}+2{a}^{3}{\it x_3}-3{
\it x_2}{a}^{2}+4a+4a{\it x_1}+2a{\it x_3}-{\it x_2}=0.$$

\item[\textbullet] Bias correction term for an estimator ${\hat{a}}$:~~
$\hat{a} ( {\hat{a}}^{8}-4{\hat{a}}^{6}+6{\hat{a}}^{4}-4{\hat{a}}^{2}+1 ) /
 ( 1+2{\hat{a}}^{2} ) ^{2}.$
\end{itemize}

\subsection{Example: log marginal model}
\label{subsec:ex_log_marginal}
Here, we consider a log marginal model. See \cite{bergsma2009} for more on marginal models.

\begin{itemize}
\item[\textbullet] Statistical model (Poisson regression):\\
$X_{ij} \displaystyle\mathop{\sim}^{ i.i.d} {\rm Po}(Np_{ij})$
%$X_{ij} \displaystyle\mathop{\sim}^{ i.i.d} {\rm Po}(N p_{ij})$
s.t. 
$p_{ij}\in (0,1)$ for $i=1,2$ and $j=1,2,3$
with model constraints:
\begin{align}
p_{11} + p_{12} + p_{13} &+p_{21} + p_{22} + p_{23} = 1,\nonumber\\
p_{11} + p_{12} + p_{13} &=p_{21} + p_{22} + p_{23},\nonumber\\
\displaystyle\frac{p_{11}/p_{21}}{p_{12}/p_{22}}&=\frac{p_{12}/p_{22}}{p_{13}/p_{23}}.\label{eq:costant_ratio}
\end{align}

Condition (\ref{eq:costant_ratio}) can appear in a statistical test of whether
acceleration of the ratio $p_{1j}/p_{2j}$ is constant.

In this case, $d=6$ and $p=3$.

\item[\textbullet] Log density w.r.t. the point mass measure on $\mathbb{Z}_{\geq 0}^6$:
$$\log f(x|p)=\log\{\prod_{ij} e^{-Np_{ij}} (Np_{ij})^{X_{ij}}\}
= -N + \sum_{ij} X_{ij} \log(Np_{ij}).$$

\item[\textbullet] The full expectation family is given by
$$
\begin{bmatrix}
X_{1} & X_{2} & X_{3} \\
X_{4} & X_{5} & X_{6} \\
\end{bmatrix}
:=
\begin{bmatrix}
X_{11} & X_{12} & X_{13} \\
X_{21} & X_{22} & X_{23} \\
\end{bmatrix}
,$$
$$
\begin{bmatrix}
\eta_{1} & \eta_{2} & \eta_{3} \\
\eta_{4} & \eta_{5} & \eta_{6} \\
\end{bmatrix}
=
N \begin{bmatrix}
p_{11} & p_{12} & p_{13} \\
p_{21} & p_{22} & p_{23} \\
\end{bmatrix}
,$$

$\theta^i=\log(\eta_i)$
%$\psi(\eta)=\sum_{i=1}^6 \eta_i$,
and
$\psi(\theta)=N$.
%$\psi(\theta)=\sum_{i=1}^6 \exp(\theta^i)$.

\item[\textbullet] The Fisher metric w.r.t. $\theta$:~~
$g_{ij}=\frac{\partial^2 \psi}{\partial\theta^i \partial\theta^j}=\delta_{ij}\eta_i$.

\item[\textbullet] Selection of the model parameters:
\begin{center}$[u_1,u_2,u_3]:=[\eta_1,\eta_3,\eta_5]$ and $[v_1,v_2,v_3]:=[\eta_2,\eta_4,\eta_6]$.\end{center}

\item[\textbullet] A set of vectors $e_i\in \mathbb{R}^6$:
\begin{center}
$e_0:=
\begin{bmatrix}
{\eta_{2}^2 (\eta_{4}-\eta_{6})}\\{-\eta_{2}^2 (\eta_{4}-\eta_{6})}\\{0}\\{-\eta_{3} \eta_{5}^2-2 \eta_{2} \eta_{4} \eta_{6}}\\{0}\\{\eta_{3} \eta_{5}^2+2 \eta_{2} \eta_{4} \eta_{6}}
\end{bmatrix}\in (\nabla_u \eta)
$,
\end{center}

$[e_1,e_2,e_3]:=
\left[
\begin{bmatrix} 
\eta_1\\ \eta_2\\ \eta_3\\ 0\\ 0 \\ 0 
\end{bmatrix},
\begin{bmatrix}
{\eta_{1} (-\eta_{1} \eta_{5}^2+\eta_{3} \eta_{5}^2)}\\{\eta_{2} (-\eta_{1} \eta_{5}^2-2 \eta_{2} \eta_{4} \eta_{6})}\\{0}\\{\eta_{4} (\eta_{2}^2 \eta_{4}-\eta_{2}^2 \eta_{6})}\\{\eta_{5} (\eta_{2}^2 \eta_{4}+2 \eta_{1} \eta_{3} \eta_{5})}\\{0}
\end{bmatrix},
\begin{bmatrix}
{\eta_{1} (\eta_{1} \eta_{5}^2-\eta_{3} \eta_{5}^2)}\\{\eta_{2} (\eta_{1} \eta_{5}^2+2 \eta_{2} \eta_{4} \eta_{6})}\\{0}\\{\eta_{4} (2 \eta_{1} \eta_{3} \eta_{5}+\eta_{2}^2 \eta_{6})}\\{0}\\{\eta_{6} (\eta_{2}^2 \eta_{4}+2 \eta_{1} \eta_{3} \eta_{5})}
\end{bmatrix}
\right]
\in ((\nabla_u \eta)^{\perp_{\bar{G}}})^3$
\vspace{0.1cm}

\item[\textbullet] A vector version of the second-order efficient estimator is, for example,
$$X-\eta+v_1\cdot e_1 + v_2 \cdot e_2 
+ v_3 \cdot e_3 + c \cdot v_1^3\cdot e_0=0.$$

\item[\textbullet] The bias correction term of the estimator = 0.

\item[\textbullet] A set of estimating equations for MLE:

\{ 
{$x_{{1}}{\eta_{{2}}}^{2}{\eta_{{4}}}^{2}\eta_{{6}}-x_{{1}}{\eta_{{2}}}^{2}\eta_{{4}}{
\eta_{{6}}}^{2}-x_{{2}}\eta_{{1}}\eta_{{2}}{\eta_{{4}}}^{2}\eta_{{6}}+x_{{2}}\eta_{{1}}\eta
_{{2}}\eta_{{4}}{\eta_{{6}}}^{2}-2\,x_{{4}}\eta_{{1}}\eta_{{2}}\eta_{{4}}{\eta_{{6}}}^{2
}-x_{{4}}\eta_{{1}}\eta_{{3}}{\eta_{{5}}}^{2}\eta_{{6}}+2\,x_{{6}}\eta_{{1}}\eta_{{2}}{\eta
_{{4}}}^{2}\eta_{{6}}+x_{{6}}\eta_{{1}}\eta_{{3}}\eta_{{4}}{\eta_{{5}}}^{2}$},\\
$-x_{{2}}\eta
_{{2}}\eta_{{3}}{\eta_{{4}}}^{2}\eta_{{6}}+x_{{2}}\eta_{{2}}\eta_{{3}}\eta_{{4}}{\eta_{{6}}
}^{2}+x_{{3}}{\eta_{{2}}}^{2}{\eta_{{4}}}^{2}\eta_{{6}}-x_{{3}}{\eta_{{2}}}^{2}\eta_{
{4}}{\eta_{{6}}}^{2}-x_{{4}}\eta_{{1}}\eta_{{3}}{\eta_{{5}}}^{2}\eta_{{6}}-2\,x_{{4}}
\eta_{{2}}\eta_{{3}}\eta_{{4}}{\eta_{{6}}}^{2}+x_{{6}}\eta_{{1}}\eta_{{3}}\eta_{{4}}{\eta_{{5}
}}^{2}+2\,x_{{6}}\eta_{{2}}\eta_{{3}}{\eta_{{4}}}^{2}\eta_{{6}}$,\\
{
$-2\,x_{{4}}\eta_{{1}}
\eta_{{3}}{\eta_{{5}}}^{2}\eta_{{6}}-x_{{4}}{\eta_{{2}}}^{2}\eta_{{4}}\eta_{{5}}\eta_{{6}}+
x_{{5}}{\eta_{{2}}}^{2}{\eta_{{4}}}^{2}\eta_{{6}}-x_{{5}}{\eta_{{2}}}^{2}\eta_{{4}}{\eta
_{{6}}}^{2}+2\,x_{{6}}\eta_{{1}}\eta_{{3}}\eta_{{4}}{\eta_{{5}}}^{2}+x_{{6}}{\eta_{{2
}}}^{2}\eta_{{4}}\eta_{{5}}\eta_{{6}}$},\\
$\eta_{{1}}\eta_{{3}}{\eta_{{5}}}^{2}-{\eta_{{2}}}^{2}
\eta_{{4}}\eta_{{6}}$,
{
$\eta_{{1}}+\eta_{{2}}+\eta_{{3}}-\eta_{{4}}-\eta_{{5}}-\eta_{{6}}$},
$-\eta_{{1}
}-\eta_{{2}}-\eta_{{3}}-\eta_{{4}}-\eta_{{5}}-\eta_{{6}}+1$\}

The total degree of the equations is $5\times 5\times 5 \times 4\times 1\times 1= 500.$

\item[\textbullet] A set of estimating equations for a 2nd-order efficient estimator with degree 
at most 2:

\begin{small}
\{{$-3\,x_{{1}}x_{{2}}{x_{{4}}}^{2}x_{{6}}\eta_{{2}}+6\,x_{{1}}x_{{2}}{x_{{4
}}}^{2}x_{{6}}\eta_{{6}}+x_{{1}}x_{{2}}{x_{{4}}}^{2}\eta_{{2}}\eta_{{6}}-2\,x_{
{1}}x_{{2}}{x_{{4}}}^{2}{\eta_{{6}}}^{2}+3\,x_{{1}}x_{{2}}x_{{4}}{x_{{6}}
}^{2}\eta_{{2}}-6\,x_{{1}}x_{{2}}x_{{4}}{x_{{6}}}^{2}\eta_{{4}}+2\,x_{{1}}x_
{{2}}x_{{4}}x_{{6}}\eta_{{2}}\eta_{{4}}-2\,x_{{1}}x_{{2}}x_{{4}}x_{{6}}\eta_{{2
}}\eta_{{6}}-x_{{1}}x_{{2}}{x_{{6}}}^{2}\eta_{{2}}\eta_{{4}}+2\,x_{{1}}x_{{2}}{
x_{{6}}}^{2}{\eta_{{4}}}^{2}+3\,x_{{1}}x_{{3}}x_{{4}}{x_{{5}}}^{2}\eta_{{6}}
-2\,x_{{1}}x_{{3}}x_{{4}}x_{{5}}\eta_{{5}}\eta_{{6}}-3\,x_{{1}}x_{{3}}{x_{{5
}}}^{2}x_{{6}}\eta_{{4}}+2\,x_{{1}}x_{{3}}x_{{5}}x_{{6}}\eta_{{4}}\eta_{{5}}+x_
{{1}}{x_{{4}}}^{2}x_{{6}}{\eta_{{2}}}^{2}-2\,x_{{1}}{x_{{4}}}^{2}x_{{6}}\eta
_{{2}}\eta_{{6}}-x_{{1}}x_{{4}}{x_{{5}}}^{2}\eta_{{3}}\eta_{{6}}-x_{{1}}x_{{4}}
{x_{{6}}}^{2}{\eta_{{2}}}^{2}+2\,x_{{1}}x_{{4}}{x_{{6}}}^{2}\eta_{{2}}\eta_{{4}
}+x_{{1}}{x_{{5}}}^{2}x_{{6}}\eta_{{3}}\eta_{{4}}+3\,{x_{{2}}}^{2}{x_{{4}}}^
{2}x_{{6}}\eta_{{1}}-{x_{{2}}}^{2}{x_{{4}}}^{2}\eta_{{1}}\eta_{{6}}-3\,{x_{{2}}
}^{2}x_{{4}}{x_{{6}}}^{2}\eta_{{1}}-2\,{x_{{2}}}^{2}x_{{4}}x_{{6}}\eta_{{1}}
\eta_{{4}}+2\,{x_{{2}}}^{2}x_{{4}}x_{{6}}\eta_{{1}}\eta_{{6}}+{x_{{2}}}^{2}{x_{
{6}}}^{2}\eta_{{1}}\eta_{{4}}-x_{{2}}{x_{{4}}}^{2}x_{{6}}\eta_{{1}}\eta_{{2}}-2\,x
_{{2}}{x_{{4}}}^{2}x_{{6}}\eta_{{1}}\eta_{{6}}+x_{{2}}x_{{4}}{x_{{6}}}^{2}\eta_
{{1}}\eta_{{2}}+2\,x_{{2}}x_{{4}}{x_{{6}}}^{2}\eta_{{1}}\eta_{{4}}-x_{{3}}x_{{4
}}{x_{{5}}}^{2}\eta_{{1}}\eta_{{6}}+x_{{3}}{x_{{5}}}^{2}x_{{6}}\eta_{{1}}\eta_{{4}
}$},\\
$3\,x_{{1}}x_{{3}}x_{{4}}{x_{{5}}}^{2}\eta_{{6}}-2\,x_{{1}}x_{{3}}x_{{4}
}x_{{5}}\eta_{{5}}\eta_{{6}}-3\,x_{{1}}x_{{3}}{x_{{5}}}^{2}x_{{6}}\eta_{{4}}+2
\,x_{{1}}x_{{3}}x_{{5}}x_{{6}}\eta_{{4}}\eta_{{5}}-x_{{1}}x_{{4}}{x_{{5}}}^{
2}\eta_{{3}}\eta_{{6}}+x_{{1}}{x_{{5}}}^{2}x_{{6}}\eta_{{3}}\eta_{{4}}+3\,{x_{{2}}
}^{2}{x_{{4}}}^{2}x_{{6}}\eta_{{3}}-{x_{{2}}}^{2}{x_{{4}}}^{2}\eta_{{3}}\eta_{{
6}}-3\,{x_{{2}}}^{2}x_{{4}}{x_{{6}}}^{2}\eta_{{3}}-2\,{x_{{2}}}^{2}x_{{4}
}x_{{6}}\eta_{{3}}\eta_{{4}}+2\,{x_{{2}}}^{2}x_{{4}}x_{{6}}\eta_{{3}}\eta_{{6}}+{x
_{{2}}}^{2}{x_{{6}}}^{2}\eta_{{3}}\eta_{{4}}-3\,x_{{2}}x_{{3}}{x_{{4}}}^{2}x
_{{6}}\eta_{{2}}+6\,x_{{2}}x_{{3}}{x_{{4}}}^{2}x_{{6}}\eta_{{6}}+x_{{2}}x_{{
3}}{x_{{4}}}^{2}\eta_{{2}}\eta_{{6}}-2\,x_{{2}}x_{{3}}{x_{{4}}}^{2}{\eta_{{6}}}
^{2}+3\,x_{{2}}x_{{3}}x_{{4}}{x_{{6}}}^{2}\eta_{{2}}-6\,x_{{2}}x_{{3}}x_{
{4}}{x_{{6}}}^{2}\eta_{{4}}+2\,x_{{2}}x_{{3}}x_{{4}}x_{{6}}\eta_{{2}}\eta_{{4}}
-2\,x_{{2}}x_{{3}}x_{{4}}x_{{6}}\eta_{{2}}\eta_{{6}}-x_{{2}}x_{{3}}{x_{{6}}}
^{2}\eta_{{2}}\eta_{{4}}+2\,x_{{2}}x_{{3}}{x_{{6}}}^{2}{\eta_{{4}}}^{2}-x_{{2}}
{x_{{4}}}^{2}x_{{6}}\eta_{{2}}\eta_{{3}}-2\,x_{{2}}{x_{{4}}}^{2}x_{{6}}\eta_{{3
}}\eta_{{6}}+x_{{2}}x_{{4}}{x_{{6}}}^{2}\eta_{{2}}\eta_{{3}}+2\,x_{{2}}x_{{4}}{
x_{{6}}}^{2}\eta_{{3}}\eta_{{4}}+x_{{3}}{x_{{4}}}^{2}x_{{6}}{\eta_{{2}}}^{2}-2
\,x_{{3}}{x_{{4}}}^{2}x_{{6}}\eta_{{2}}\eta_{{6}}-x_{{3}}x_{{4}}{x_{{5}}}^{2
}\eta_{{1}}\eta_{{6}}-x_{{3}}x_{{4}}{x_{{6}}}^{2}{\eta_{{2}}}^{2}+2\,x_{{3}}x_{
{4}}{x_{{6}}}^{2}\eta_{{2}}\eta_{{4}}+x_{{3}}{x_{{5}}}^{2}x_{{6}}\eta_{{1}}\eta_{{
4}}$,\\
{
$6\,x_{{1}}x_{{3}}x_{{4}}{x_{{5}}}^{2}\eta_{{6}}-4\,x_{{1}}x_{{3}}x_{{
4}}x_{{5}}\eta_{{5}}\eta_{{6}}-6\,x_{{1}}x_{{3}}{x_{{5}}}^{2}x_{{6}}\eta_{{4}}+
4\,x_{{1}}x_{{3}}x_{{5}}x_{{6}}\eta_{{4}}\eta_{{5}}-2\,x_{{1}}x_{{4}}{x_{{5}
}}^{2}\eta_{{3}}\eta_{{6}}+2\,x_{{1}}{x_{{5}}}^{2}x_{{6}}\eta_{{3}}\eta_{{4}}+3\,{
x_{{2}}}^{2}{x_{{4}}}^{2}x_{{6}}\eta_{{5}}-{x_{{2}}}^{2}{x_{{4}}}^{2}\eta_{{
5}}\eta_{{6}}-3\,{x_{{2}}}^{2}x_{{4}}x_{{5}}x_{{6}}\eta_{{4}}+3\,{x_{{2}}}^{
2}x_{{4}}x_{{5}}x_{{6}}\eta_{{6}}+{x_{{2}}}^{2}x_{{4}}x_{{5}}\eta_{{4}}\eta_{{6
}}-{x_{{2}}}^{2}x_{{4}}x_{{5}}{\eta_{{6}}}^{2}-3\,{x_{{2}}}^{2}x_{{4}}{x_
{{6}}}^{2}\eta_{{5}}-{x_{{2}}}^{2}x_{{4}}x_{{6}}\eta_{{4}}\eta_{{5}}+{x_{{2}}}^
{2}x_{{4}}x_{{6}}\eta_{{5}}\eta_{{6}}+{x_{{2}}}^{2}x_{{5}}x_{{6}}{\eta_{{4}}}^{
2}-{x_{{2}}}^{2}x_{{5}}x_{{6}}\eta_{{4}}\eta_{{6}}+{x_{{2}}}^{2}{x_{{6}}}^{2
}\eta_{{4}}\eta_{{5}}-2\,x_{{2}}{x_{{4}}}^{2}x_{{6}}\eta_{{2}}\eta_{{5}}+2\,x_{{2}
}x_{{4}}x_{{5}}x_{{6}}\eta_{{2}}\eta_{{4}}-2\,x_{{2}}x_{{4}}x_{{5}}x_{{6}}\eta_
{{2}}\eta_{{6}}+2\,x_{{2}}x_{{4}}{x_{{6}}}^{2}\eta_{{2}}\eta_{{5}}-2\,x_{{3}}x_
{{4}}{x_{{5}}}^{2}\eta_{{1}}\eta_{{6}}+2\,x_{{3}}{x_{{5}}}^{2}x_{{6}}\eta_{{1}}
\eta_{{4}}$,
}\\
$\eta_{{1}}\eta_{{3}}{\eta_{{5}}}^{2}-{\eta_{{2}}}^{2}\eta_{{4}}\eta_{{6}}$,
{
$\eta_{{1}}+\eta_{{2}}+\eta_{{3}}-\eta_{{4}}-\eta_{{5}}-\eta_{{6}}$
},
$-\eta_{{1}}-\eta_{{2}}-\eta_{{3}}-\eta_{{4}}-\eta_{{5}}-\eta_{{6}}+1$\}
\end{small}
The total degree of the polynomial equations is $32$

\item[\textbullet] A set of estimating equations for a first-order-efficient estimator with degree 
at most 1:

%\begin{small}
{\{$-{x_{{5}}}^{2}x_{{4}}\eta_{{6}}x_{{1}}x_{{3}}+{x_{{5}}}^{2}x_{{6}}\eta_{{4}
}x_{{1}}x_{{3}}+2\,{x_{{6}}}^{2}\eta_{{4}}x_{{1}}x_{{2}}x_{{4}}-2\,{x_{{4
}}}^{2}\eta_{{6}}x_{{1}}x_{{2}}x_{{6}}-{x_{{6}}}^{2}x_{{1}}x_{{2}}\eta_{{2}}
x_{{4}}+{x_{{4}}}^{2}x_{{1}}x_{{2}}\eta_{{2}}x_{{6}}+{x_{{2}}}^{2}{x_{{6}
}}^{2}\eta_{{1}}x_{{4}}-{x_{{4}}}^{2}{x_{{2}}}^{2}\eta_{{1}}x_{{6}}
$
},\\
$
-{x_{{5}
}}^{2}x_{{4}}\eta_{{6}}x_{{1}}x_{{3}}+{x_{{5}}}^{2}x_{{6}}\eta_{{4}}x_{{1}}x
_{{3}}+2\,{x_{{6}}}^{2}\eta_{{4}}x_{{2}}x_{{3}}x_{{4}}-2\,{x_{{4}}}^{2}\eta_
{{6}}x_{{2}}x_{{3}}x_{{6}}-{x_{{6}}}^{2}x_{{2}}x_{{3}}\eta_{{2}}x_{{4}}+{
x_{{4}}}^{2}x_{{2}}x_{{3}}\eta_{{2}}x_{{6}}+{x_{{2}}}^{2}{x_{{6}}}^{2}\eta_{
{3}}x_{{4}}-{x_{{4}}}^{2}{x_{{2}}}^{2}\eta_{{3}}x_{{6}}
$,\\
{
$
-2\,{x_{{5}}}^{2}
x_{{4}}\eta_{{6}}x_{{1}}x_{{3}}+2\,{x_{{5}}}^{2}x_{{6}}\eta_{{4}}x_{{1}}x_{{
3}}-x_{{4}}x_{{6}}x_{{5}}{x_{{2}}}^{2}\eta_{{6}}+x_{{4}}x_{{5}}{x_{{2}}}^
{2}\eta_{{4}}x_{{6}}-{x_{{4}}}^{2}{x_{{2}}}^{2}\eta_{{5}}x_{{6}}+x_{{4}}{x_{
{6}}}^{2}{x_{{2}}}^{2}\eta_{{5}}
$
},\\
$\eta_{{1}}\eta_{{3}}{\eta_{{5}}}^{2}-{\eta_{{2}}}^{2}\eta_{{4}}\eta_{{6}}$,
{
$\eta_{{1}}+\eta_{{2}}+\eta_{{3}}-\eta_{{4}}-\eta_{{5}}-\eta_{{6}}$
},
$-\eta_{{1}}-\eta_{{2}}-\eta_{{3}}-\eta_{{4}}-\eta_{{5}}-\eta_{{6}}+6$\}

%\end{small}
\end{itemize}

The estimating equations for a 2nd-order-efficient estimator above
look much more complicated than the estimating equation for the MLE, 
but each term of the first three polynomials are at most degree 2.
Thanks to this degree reduction, the computational costs for
the estimates become much smaller as we will see in the next section.

\section{Computation}
To obtain estimates based on the method of this paper,
we need fast algorithms to find the solution of polynomial equations.
The authors have carried out computations using homotopy continuation method (matlab program HOM4PS2
by Lee, Li and Tsuai \cite{lee_etal2008})
for the log marginal model in Sec. \ref{subsec:ex_log_marginal} and
a data $\bar{X}=(1,1,1,1,1,1)$. 

The run time to compute each estimate on a standard laptop 
(Intel(R) Core (TM) i7-2670QM CPU, 2.20GHz, 4.00GB memory) 
is given by Table \ref{tb:homotopy-time}.
The computation is repeated 10 times and the averages and the standard deviations are displayed.
Note the increasing of the speed for the second-order efficient estimators
is due to the degree reduction technique.
The term ``path'' in the table heading refers to a primitive iteration step
within the homotopy method.
In the faster polyhedron version, the solution region is subdivided into polyhedral domains.

% {Computational Results by
% the Homotopy Continuation Methods}
% \begin{itemize}
% \item
% Software for the homotopy methods: HOM4PS2 by Lee, Li and Tsuai.
% \item
% $X=(1,1,1,1,1,1).$
% \item
% Repeat count: 10.
% \end{itemize}
\begin{table}
\begin{center}
\caption{Computational time for each estimate by the homotopy continuation methods}
\begin{tabular}{|l|l|r|c|}
\hline
 algorithm   & estimator &   \#paths & running time [s]\\ 
 &&&(avg. $\pm$ std.) \\ \hline
 linear        & MLE        & 500         & 1.137 $\pm$ 0.073 \\ \cline{2-4}
 homotopy  & 2nd eff.    & 32           & {\bf 0.150   $\pm$ 0.047}  \\ \hline
 polyhedral & MLE         & 64           & 0.267 $\pm$ 0.035  \\ \cline{2-4}
  homotopy &2nd eff     & 24           & {\bf 0.119 $\pm$ 0.027}  \\ \hline
\end{tabular}
\label{tb:homotopy-time}
\end{center}
\end{table}

\begin{figure}[htbp]
\begin{center}
 \begin{minipage}{6cm}
  \begin{center}
   \includegraphics[width=6cm]{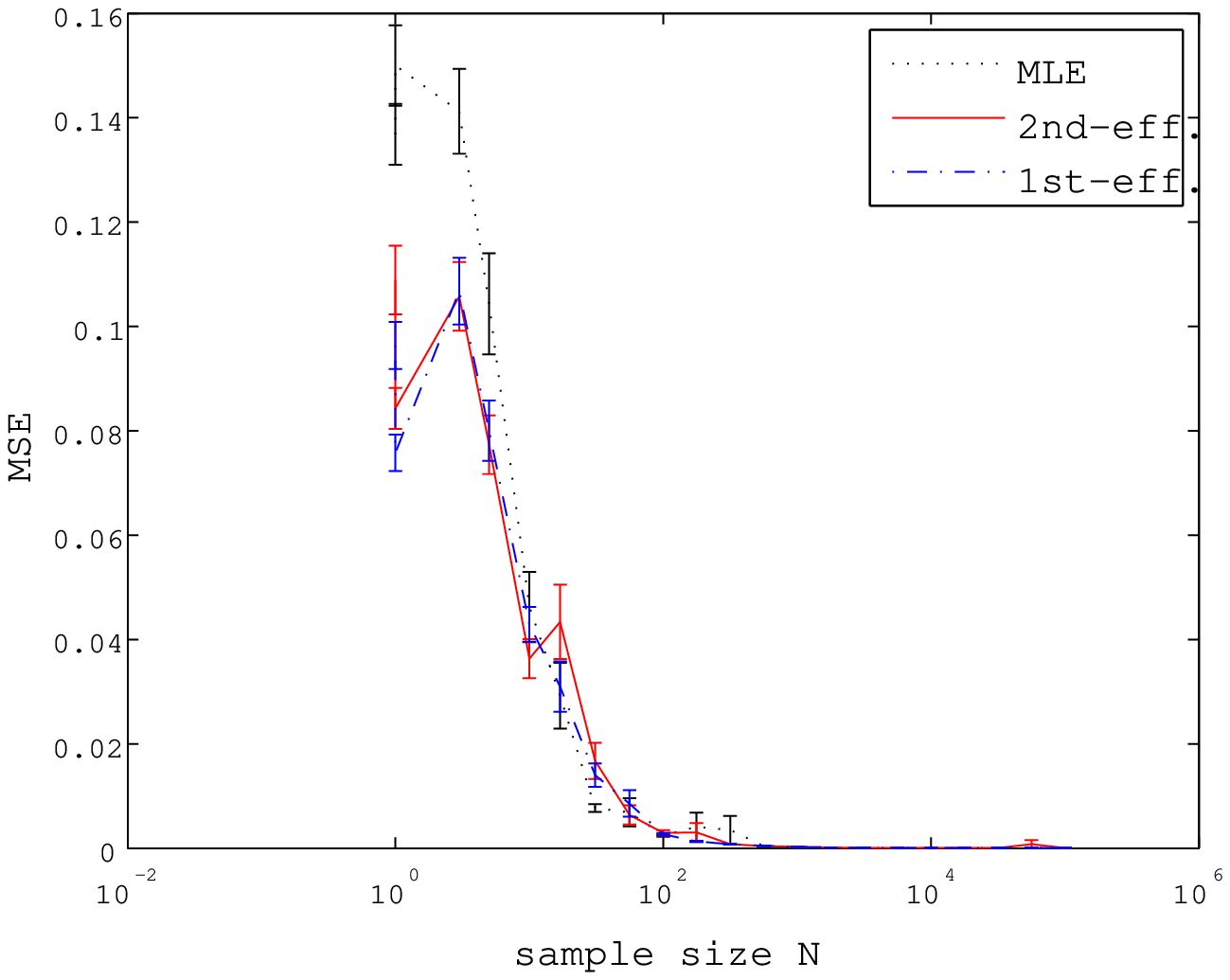}
   (1) Mean squared error
  \end{center}
  %\caption{Mean squared error}
  %\label{fig:homotopy-err}
 \end{minipage}
 \begin{minipage}{6cm}
  \begin{center}
   \includegraphics[width=6cm]{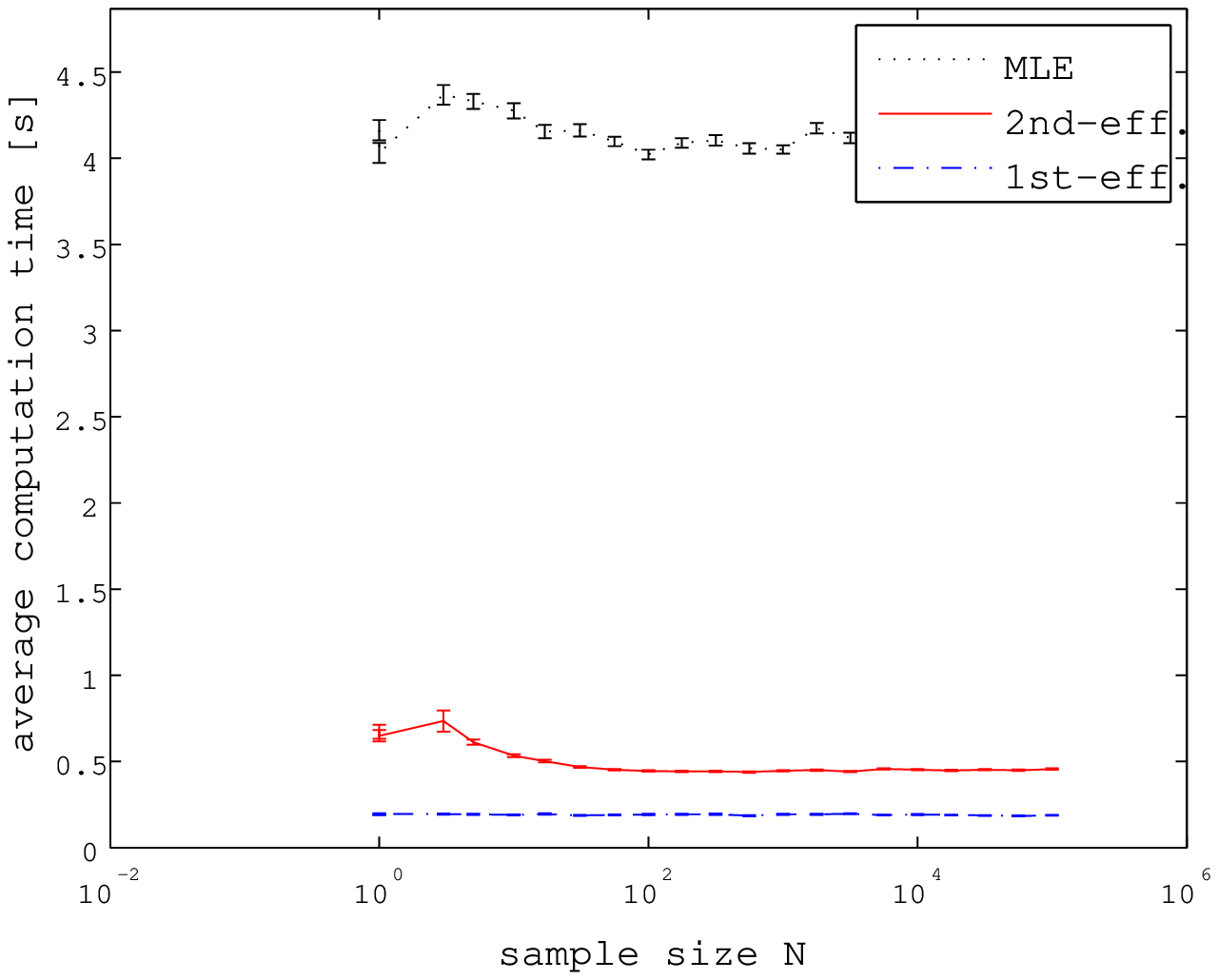}
   (2) Computation time (s)
  \end{center}
  %\caption{Computation time (s)}
  %\label{fig:homotopy-time}
 \end{minipage}
\caption{The mean squared error and computation time for each estimate by the homotopy continuation method} 
\label{fig:homotopy-err-time}
\end{center}
\end{figure}

Figure \ref{fig:homotopy-err-time} shows the mean squared error and
the computational time of the MLE, the first-order estimator and the second-order efficient estimator
of Sec.~\ref{subsec:ex_log_marginal}.
The true parameter is set $\eta^*=(1/6,1/4,1/12,1/12,1/4,1/6)$, a point in the model manifold,
and $N$ random samples are generated i.i.d. from the distribution with the parameter.
The computation is repeated for exponentially increasing sample sizes $N=1,...,10^5$.
%More extensive evaluation will appear in the longer version of this paper
In general, there are multiple roots for polynomial equations and here
we selected the root closest to the sample mean by the Euclidean norm.
Figure \ref{fig:homotopy-err-time} (1) also shows that the mean squared error is approximately
the same for the three estimators, but (2) shows that the computational time is much more for
the MLE. 

\section{Discussion}
In this paper we have concentrated on reduction of the polynomial degree of the estimating equations and shown the benefits in computation of the solutions. We do not expect the estimators to be closed form, such as a rational polynomial form in the data. The most we can expect is they are algebraic, that is they are the solution of algebraic equations. They lie on a zero dimensional algebraic variety. It is clear that there is no escape from using mixed symbolic-numerical methods. In algebraic statistics the number of solution of the ML equation is called the ML degree. Given that we have more general estimating equations than pure ML equations this points to an extended theory or ``quasi'' ML degree of efficient estimator degree.  The issue of exactly which solution to use as our estimators persists. In the paper we suggest taking the solution closest to the sufficient statistic in the Euclidian metric. We could use other metrics and more theory is needed.

Here we have put forward estimating equations with reduced degree and shown the benefits in terms of computation. But we could have used other criteria for choosing the equations, while remaining in the efficient class. We might prefer to choose an equation which reduces the bias further via decreasing the next order term. There may thus be some trade off between degree and bias.

Beyond the limited ambitions of this paper to look at second-order efficiency lie several other areas, notably
hypothesis testing and model selection. But the question is the same: to what extent can we bring the algebraic methods to bear, for example by expressing additional differential forms and curvatures in algebraic
terms. Although estimation typically requires a mixture of symbolic and numeric methods  in some cases only the computation of the efficient estimate requires numeric procedures and the other computations can be carrying out symbolically.

\appendix

% \section{Asymptotic estimation theory via information geometry (Summary)}
% We give a short summary of asymptotic estimation theory in \cite{amari1985}
% especially on curved exponential families.
% Note that in order to verify the following asymptotic theory,
% we assume some regularity conditions on statistical distribution families
% to guarantee existence of the moments $E_\theta[X^p]$ with up to some sufficient degree
% and swapping of the integration and the partial derivative as $\frac{\partial}{\partial \theta} E_\theta[g(x,\theta)]
% =E_\theta[\frac{\partial}{\partial \theta}g(x,\theta)]$ for any function $g$ appeared in
% asymptotic expansions. See section 2.2 of  \cite{amari1985} for more details.

% \begin{theorem}[Theorem 5.4 of \cite{amari1985}], 
% \label{thm:amari-1st-eff}
% $E_u[(\hat{u}^a-u^a)(\hat{u}^b-u^b)]=N^{-1}[g_{ab}-g_{a\kappa}g^{\kappa \lambda}g_{b\lambda}]^{-1}+O(N^{-2})$.\\
% Thus, an estimator is \begin{it} 1-st order efficient \end{it}
% iff $g_{a\kappa}=0$.
% \end{theorem}

\section{Normal forms}
A basic text for the materials in this section is \cite{cox2007}.
The rapid growth of modern computational algebra can be credited
to the celebrated Buchberger's algorithm \cite{buchberger2006}.

A monomial ideal $I$ in a polynomial ring $K[x_1,\dots,x_n]$
over a field $K$ is an ideal for which there is a collection of
monomials $f_1,\ldots, f_m$ such that any $g \in I$ can be expressed
as a sum $$g = \sum_{i=1}^m g_i(x)f_i(x)$$
with some polynomials $g_i\in K[x_1,\dots,x_n]$.
We can appeal to the representation of a monomial $x^{\alpha}=x_1^{\alpha_1}\dots
x_n^{\alpha_n}$ by
its exponent $\alpha=(\alpha_1,\dots,\alpha_n)$. 
If $\beta \geq 0$ is another exponent then
$$x^{\alpha} x^{\beta} = x^{\alpha + \beta},$$
and $\alpha + \beta$ is in the positive (shorthand for non-negative)
``orthant" with corner at $\alpha$. The set of all monomials in a
monomial ideal is the union of all positive orthants whose ``corners"
are given  by the exponent vectors of  the generating monomial $f_1,
\ldots, f_m$.  A monomial ordering written $x^{\alpha} \prec x^{\beta}$ is a total (linear) ordering
on monomials such that for $\gamma \geq 0$,  
$x^{\alpha} \prec x^{\beta} \Rightarrow  x^{\alpha +\gamma} \prec x^{\beta+ \gamma}$. 
Any polynomial $f(x)$ has a leading terms with respect to $\prec$, written $LT(f)$.

There are, in general,  many ways to express a given ideal $I$ as
being generated from a basis $I = \langle  f_1,\dots, f_m\rangle$.
That is to say, there are many choices of basis. Given an ideal $I$ a set $\{g_1, \ldots g_m\}$ is called a Gr\"obner
basis (G-basis) if:
$$\langle LT(g_1), \ldots, LT(g_m)\rangle\; = \; \langle LT(I)\rangle,$$
where $\langle LT(I)\rangle$ is the ideal generated by all the
monomials in $I$. We sometimes refer to $\langle LT(I) \rangle$ as the {\em leading term ideal}. 
Any ideal $I$ has a  Gr\"obner basis and any Gr\"obner basis in the ideal is a basis of the ideal.

Given a monomial
ordering and an ideal expressed in terms of the G-basis, $I\;=\;
\langle  g_1,\ldots, g_m\rangle$, any polynomial $f$ has a unique
remainder with respect the quotient operation $K[x_1, \ldots,
x_k]/I$. That is
$$f = \sum_{i=1}^m s_i(x)g_i(x) + r(x).$$
We call the remainder $r(x)$ the {\em normal form} of $f$ with
respect to $I$ and write $NF(f)$. Or, to stress the fact that it may
depend on $\prec$, we write $NF(f, \prec)$. Given a monomial ordering $\prec$, a polynomial $f=\sum_{\alpha \in
L} \theta_{\alpha} x^{\alpha}$ for some $L$ is a normal form with respect to
$\prec$ if $x^{\alpha} \notin \langle LT(f) \rangle $ for all
$\alpha \in L$. An equivalent way of saying this is: given an ideal $I$  and a monomial ordering $\prec$, for every $f
\in K[x_1,\ldots,x_k]$ there is a unique normal form $NF(f)$ such
that $f-NF(f) \in I$.

\section{Homotopy continuation method}
Homotopy continuation method is an algorithm to find the solutions of simultaneous polynomial equations numerically. 
See, for example, \cite{verschelde1999} and \cite{li1997} for more details of the algorithm and theory.

We will explain the method briefly by a simple example of 2 equations with 2 unknowns
\begin{center}
Input: $f, g \in \mathbb{R}[x,y]$
\end{center}
\begin{center}
Output: The solutions of $f(x,y)=g(x,y)=0$
\end{center}

\begin{figure}
\begin{center}
\includegraphics[width=6.5cm]{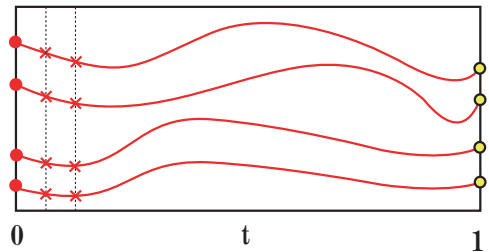}
\caption{Paths for the homotopy continuation method}
\label{fig:homotopy-2}
\end{center}
\end{figure}

\begin{itemize}
\item[]
{\bf Step 1} Select arbitrary polynomials of the form:
\begin{align}
f_0(x,y):=f_0(x):=a_1x^{d_1}-b_1&=0,\nonumber\\
g_0(x,y):=g_0(y):=a_2y^{d_2}-b_2&=0\label{eq:homotopy}
\end{align}
where $d_1= \deg(f)$ and $d_2=\deg(g)$.
Polynomial equations in this form are easy to solve.

\item[]
{\bf Step 2} Take the convex combinations:
\begin{align*}
f_t(x,y):=t f(x,y)+ (1-t) f_0(x,y),\\
g_t(x,y):=t g(x,y)+ (1-t) g_0(x,y)
\end{align*}
then our target becomes the solution for $t=1$.

\item[]
{\bf Step 3} 
Compute the solution for $t=\delta$ for small $\delta$ by the solution for $t=0$ numerically.

\item[]
{\bf Step 4}
Repeat this until we obtain the solution for $t=1$.\\
\end{itemize}

Figure \ref{fig:homotopy-2} shows a sketch of the algorithm.
This algorithm is called {\em the (linear) homotopy continuation method}
and justified if the path connects $t=0$ and $t=1$ continuously
without an intersection. That can be proved for almost all $a$ and $b$. See \cite{li1997}.

For each computation for the homotopy continuation method, 
the number of the paths is the number of the solutions of
(\ref{eq:homotopy}).
%\begin{align*}
%f_0(x,y):=f_0(x):=a_1x^{d_1}-b_1&=0,\\
%g_0(x,y):=g_0(y):=a_2y^{d_2}-b_2&=0.
%\end{align*}
In this case, the number of paths is $d_1 d_2$.
In general case with $m$ unknowns, it becomes $\prod_{i=1}^m d_i$ and
this causes a serious problem for computational cost.
Therefore decreasing the degree of second-order efficient estimators plays an important role for the homotopy continuation method.

Note that in order to solve this computational problem, the authors of \cite{huber1995}
proposed
the nonlinear homotopy continuation methods (or the polyhedral continuation methods).
But as we can see in Section \ref{subsec:ex_log_marginal},
the degree of the polynomials still affects the computational costs.

\section*{Acknowledgment}
This paper has benefited from conversations with and advice from a number of colleagues.
We should thank Satoshi Kuriki, Tomonari Sei, Wicher Bergsma and Wilfred Kendall.
The first author acknowledges support by JSPS KAKENHI Grant 20700258, 24700288 and
the second author acknowledges support from the Institute of Statistical Mathematics
for two visits in 2012 and 2013 and from UK EPSRC Grant EP/H007377/1.
A first version of this paper was delivered at the WOGAS3 meeting at the University of Warwick
in 2011. We thank the sponsors.
The authors also thank the referees of the short version in GSI2013 and the referees of
the first long version of the paper for insightful suggestions.

\thispagestyle{empty}


\begin{thebibliography}{99}
\addcontentsline{toc}{chapter}{References}

\bibitem{adler}
Adler, R.J. and Taylor, J.E.,
``Random Fields and Geometry'', Series: Springer Monographs in Mathematics, Springer, 2007

\bibitem{amari1982}
Amari, S.,
``Differential geometry of curved exponential families-curvatures and information loss'',
{\it Ann. Statist.}, pp. 357-385, 1982.

\bibitem{amari1985}
Amari, S., ``Differential-geometrical methods in statistics,'' Springer, 1985.

\bibitem{amari1983}
Amari, S. and Kumon, M.,
``Differential geometry of Edgeworth expansions in curved exponential family'',
{\it Ann. Inst. Stat. Math}. vol.~35, no.~1, pp. 1-24, 1983.

\bibitem{amari-nagaoka2007}
Amari, S. and Nagaoka, H., ``Methods of Information Geometry,'' vol. 191, Amer Mathematical Society, 2007.

\bibitem{bergsma2009}
Bergsma, W.P., Croon, M. and Hagenaars, J.A.,
`` Marginal Models for Dependent, Clustered, and Longitudinal Categorical Data'',
Springer,  2009.

\bibitem{buchberger2006}
Buchberger, B.,
``Bruno Buchbergerfs PhD thesis 1965: An algorithm for finding the basis elements of the residue class ring of a zero dimensional polynomial ideal'', 
{\it Journal of symbolic computation}, vol.~41, no.~3, pp. 475-511, 2006.

\bibitem{cox2007}
Cox, D.A., Little, J. and O'Shea, D.,
`` Ideals, Varieties, and Algorithms: An Introduction to Computational Algebraic Geometry and Commutative Algebra, 3/e
(Undergraduate Texts in Mathematics)'', Springer-Verlag New York, Inc., 2007.

\bibitem{dawid1977}
Dawid, A.P.,
``Further comments on some comments on a paper by Bradley Efron''. 
{\it Ann. Statist.},.vol.~5, no.~6 p.1249, 1977.

\bibitem{drton}
Drton, M.,
``Likelihood ratio tests and singularities'',
{\it Ann. Statist.}, pp. 979-1012, 2009.

\bibitem{drton_book}
Drton, M., Sturmfels B. and Sullivant, S.,
``Lectures on algebraic statistics'',
Springer, 2009.

\bibitem{efron1975}
Efron, B.,
``Defining the curvature of a statistical problem 
(with applications to second order efficiency)'',
{\it Ann. Statist.}, vol.~3, pp. 1189-1242, 1975.

\bibitem{gehrmann-2012}
Gehrmann, H., and Lauritzen. S.L.,
``Estimation of means in graphical Gaussian models with symmetries.''
{\it Ann. Statist.}, vol.~40, no.~2, pp. 1061-1073, 2012.

\bibitem{gibilisco_etal2009}
Gibilisco, P., Riccomagno, E., Rogantin, M.P. and Wynn, H.P., 
``Algebraic and geometric methods in statistics'', Cambridge University Press, 2009

\bibitem{huber1995}
Huber, B. and Sturmfels, B.,
``A polyhedral method for solving sparse polynomial systems'',
{\it Math. Comput.}, vol.~64, pp. 1541-1555, 1995.

\bibitem{kuriki2002}
Kuriki, S. and Takemura, A.,
``On the equivalence of the tube and Euler characteristic methods for the distribution of the maximum of Gaussian fields over piecewise smooth domains'',
{\it Ann. Appl. Probab.}, vol.~12, no.~2 
pp. 768-796, 2002.

\bibitem{lee_etal2008}
Lee, T.L., Li, T.Y. and Tsai, C.H.,
``HOM4PS2.0: a software package for solving polynomial systems by the polyhedral homotopy continuation method,''
{\it Computing}, vol.~83, no.~2-3, pp. 109-133, 2008

\bibitem{li1997}
Li, T.Y.,
``Numerical solution of multivariate polynomial systems by homotopy continuation methods'',
{\it Acta numerica}, vol.~6, no.~1, pp. 399-436, 1997.

\bibitem{naiman}
Naiman, D.Q.,
``Conservative Confidence Bands in Curvilinear Regression'',
{\it Ann. Statist.}, vol.~14, no.~3, pp. 896-906, 1986.

\bibitem{pistone1995}
Pistone, G., and Sempi, C.,
``An infinite-dimensional geometric structure on the space of all the probability measures equivalent to a given one'',
{\it Ann. Statist.}, pp. 1543-1561, 1995.

\bibitem{pistone1996}
Pistone, G. and Wynn, H.P.,
``Generalised confounding with {G}r\"obner bases'',
{\it Biometrika}, vol.~83, pp. 653-666, 1996.

\bibitem{rao1945}
Rao, R.C. ``Information and accuracy attainable in the estimation of statistical parameters'', 
{\it Bulletin of the Calcutta Mathematical Society}, vol.~37, no.~3, pp. 81-91, 1945.

\bibitem{andersson-1998}
Andersson, S., and Madsen, J.,
``Symmetry and lattice conditional independence in a multivariate normal distribution."
{\it Ann. Statist.}, vol.~26, no.~2, pp. 525-572, 1998.

\bibitem{verschelde1999}
Verschelde, J.,
``Algorithm 795: PHCpack: A general-purpose solver for polynomial systems by homotopy continuation'',
{\it ACM Transactions on Mathematical Software (TOMS)}
, vol.~25, no.~2, pp. 251-276, 1999.

\bibitem{watanabe}
Watanabe, S.,
``Algebraic analysis for singular statistical estimation'',{\it Algorithmic Learning Theory}, Springer Berlin Heidelberg, 1999.

\bibitem{weyl1939}
Weyl, H.,
``On the {V}olume of {T}ubes'',
{\it Amer. J. Math.},
vol.~61, pp. 461-472, 1939.

\end{thebibliography}
\end{document}